\documentclass [a4paper,14pt ]{extarticle}
\usepackage[dvips,pdftex]{graphicx,graphics}
\usepackage{amsfonts}
\usepackage[utf8]{inputenc}
\usepackage{graphicx}
\begin{document}
 
\baselineskip=16pt
\clearpage

\newcommand{\N}{\mathbb{N}}
\newcommand{\qq}{\qquad}
\newcommand{\q}{\quad}

\newcommand{\si}{\sigma}
\newcommand{\la}{\lambda}
\newcommand{\al}{\alpha}
\newcommand{\be}{{\beta}}
\newcommand{\ve}{{\varepsilon}}
\newcommand{\vk}{{\varkappa}}

\baselineskip=16pt

 \setcounter{page}{1}

{\ }

\vspace{3ex} 

 MSC-2010-class:  11N37 (Primary); 11M26 (Secondary)

\begin{center}

\vspace{5ex}

{\ } {\bf  \Large RH-Dependent Estimates of Remainder\\
\vspace{1ex}
{\ } {\bf  \Large in Modified Mertens Formula}\footnote{\q This work
 was partly supported
by the grant of Russian Foundation of Fundamental Research\\
$\qq {\ } {\qq } {\ } $ (project \# $14-01-00684.$)}$^)$
}

\vspace{7ex}

{\bf \large by Gennadiy Kalyabin\footnote{\q Samara, 
 Russia; \ gennadiy.kalyabin@gmail.com}$^)$}

\end{center}
\vspace{3ex} 

{\it Abstract:}  Assuming the validity of Riemann Hypothesis {\bf(RH)}, 
we derive the explicit  bilateral estimates  ("narrow passage") for the
remainder in the {\it modified} Mertens asymptotic  formula
for the sums of primes' reciprocals.

These results are reversable, thus yielding some new criteria for {\bf RH}.

\vspace{1ex}


\vspace{1ex}

{\it Bibliography}: 7 items

\vspace{3ex}

{\bf 1. \ Notations, brief history  and main results}

\vspace{2ex}

 As usually, let $N, j, k, m,  n$ (perhaps with indices) run the set $\N$ of all positive integers,  $\N_0:=\N \cup\{0 \}$, $p$ run the set 
  $\mathbb{P}:=\{p_1, p_2, \dots \},\ p_j<p_{j+1},$ of all primes, 
$\ve$ be an arbitrary positive number, $\delta_k$ denote  sequences,
which $\to+0$ (perhaps different even within one and the same formula);
$\ C_(a)$ stand for positive constants which may depend only on a parameter $a$; \ symbols $\triangleright$ and $\Box$ denote the proof's beginning and end; \  $\log x$ and $\gamma$ stand (resp.) for the natural logarithm of a positive $x$ and the Euler-Masceroni constant:
  $$ \gamma:= \lim_{n\to\infty} \left( \sum_{k=1}^n \frac{1}{k}
   - \log n\right)=0.577\ 215 \dots
  \eqno{(1.1)}
  $$
\vspace{-1ex}  

In 1874 F. Mertens [1] proved his famous asymptotic formula
 $$ S(x):=\sum_{p\le x} \log \frac{p}{p-1} =\log \log x +\gamma +R(x)
\ \hbox{with } R(x)=O\left(\frac{1}{\log x}\right).
  \eqno{(1.2)}
  $$

\vspace*{\fill}
\clearpage

In 1984 assuming {\bf RH} G. Robin [2, Th. 3] has come to the fundamentally stronger estimate:
$$ |R(x)|\le \frac{\log x}{8\pi\sqrt{x}},   \q x> X_0.
  \eqno{(1.3)}
  $$

J.-L. Nicolas (1983) has considered the modified Mertens
 formula\footnote{\q For $x<3 \ \log\log\theta(x)$ cannot be defined
as a real number}$^)$:
$$ S(x)=\log\log \theta(x) +\gamma +Q(x), \q x\ge3.\q 
 \eqno{(1.4)}
  $$
which differs from (1.2) by replacing $x$ in $\log\log$ by the first
Chebyshev function: $\theta(x):=\sum \{\log p: \ p\le x\}$, 
cf. [5, S.\, 3.3], and the following assertion was established [3, th. 3]:

\vspace{1ex}

{\bf Proposition 1.} {\it \ Each of the two conditions is {\bf sufficient}
for {\bf RH}}:
$$ {\bf (i)} \q \forall \ve>0:\q B_{\ve}^+:=
 \limsup\, Q(x)\, x^{0.5-\ve}<+\infty,
$$
$$ {\bf (ii)} \q \forall \ve>0:\q B_{\ve}^-:=
 \liminf\, Q(x)\, x^{0.5-\ve}>-\infty.
 \eqno{(1.5)}
  $$

\vspace{1ex}

Later the author, basing on the connection between (1.4) and 
 the Ramanujan inequality for Gronwall numbers, has proved in [4]:

\vspace{1ex}

{\bf Proposition 2.} {\it The relationship
$$ A^+:= \limsup \, Q(x)\, \sqrt{x} \, \log x <\infty,
\eqno{(1.6)}
  $$
is necessary for {\bf RH}, and this being the case, 
then necessarily} $A^+\le2\sqrt{2}$.

\vspace{1ex}

The aim of this paper is to strengthen these results  (using different approach)   by conditional ({\bf RH}{\it-dependent}) narrow estimates for $A^+$ and for the quantity:
$$ A^-:= \liminf \, Q(x)\, \sqrt{x} \, \log x.
\eqno{(1.7)}
  $$

\vspace{1ex}

{\bf Theorem.} {\it \q Assume {\bf RH}. Then one has:}
$$ {\bf (i)}\q A^+\le 2.5, \qq {\bf (ii)}\q A^-\ge 1.5.
\eqno{(1.8)}
  $$

\vspace{1ex}

{\bf Remark 1.}\ Directly from definitions, Propositions 1, 2
 and the Theorem  it follows that 
$$\hspace{-1ex}
{\bf RH}\Rightarrow(1.8){\bf(ii)}  \Rightarrow (1.5){\bf(ii)}
 \Rightarrow {\bf RH} \Rightarrow (1.8){\bf(i)} \Rightarrow (1.6) 
\Rightarrow (1.5){\bf(i)} \Rightarrow {\bf RH}.
  $$

Therefore the Theorem demonstrates, that 
all six relationships in this chain  are equivalent.

\vspace{1ex}

Some other  criteria for {\bf RH} in diverse terms
 will be adduced in Sect. 3.

\vspace*{\fill}
\clearpage

{\bf 2. Proof of the Theorem}

\vspace{2ex}

We begin with some preliminary assertions.

 Recall that {\bf RH}
 may be reformulated in terms of the deviation of Chebyshev 
function $\theta(x)$  from $x$:
$\Delta(x):= \theta(x) - x$, which  is
{\it unconditionally} $ O(x \exp(- c\sqrt{\log x})$
by virtue of  {\bf PNT}, cf [5], th. 5.19. 

Therefore (cf. also (1.2), (1.4)) for some $\tilde{x}$
 in between of $x$ and $\theta(x)$ one has:
\vspace{-2ex}
$$ Q(x)=R(x) +\log\log x - \log\log \theta(x)
 = R(x) - \frac{\Delta(x)}{\tilde{x}\log\tilde{x}}
=O\left(\frac{1}{\log x}\right).
\eqno{(2.1)}
  $$

The  proposition  below is  well known classical result by H. von Koch (1901), cf., e.g.,  [5], th. 5.21.

\vspace{2ex}

{\bf Proposition 3.}\ {\it The following assertions are equivalent:}
$$ {\bf (i) \q RH} \hbox{ holds true}; \q  {\bf (ii)}\q 
 |\Delta(x)| <x^{0.5+\ve},\ \forall\ve>0, \ \forall x>X_{\ve};
$$
$$   {\bf (iii)} \q |\Delta(x)|\le \frac{\sqrt{x}\log^2x}{8\pi}, \ x>X_0.
\eqno{(2.2)}
  $$

\vspace{2ex}

In the Ingham’s monograph [6], S.\ V.10, form. (35), the helpful
conditional estimate for the primitive of $(\psi(x) - x)$ is adduced.
\vspace{2ex}

{\bf Proposition 4}.\ {\it Assume {\bf RH}; then}
$\left|\int_0^x (\psi(t) - t)\,dt\right| < 0.1x^{3/2}, \ x<X_0$.

\vspace{2ex}

Whence, taking into account the relationship 
$\theta(x)=\psi(x)-\sqrt{x} +O(x^{1/3})$, one obtains:


\vspace{2ex}

{\bf Proposition 5}.\ {\it The validity of {\bf RH} implies the following bilateral estimates for the $(\theta(x) - x)$-primitive}:
\vspace{-1ex}
$$ \Phi(x):=\int_0^x \Delta(t)\,dt = (-2/3 + b(x)) x^{3/2};
 \q |b(x)|<0.1, \ x>X_0.
\eqno{(2.3)}
  $$

The main point of the Theorem’s proof is the following 
corollary of {\bf RH}.

\vspace{2ex}

{\bf Lemma 1.} \ {\it Assume {\bf RH}. Introduce the quantities }
 $(3\le x< y<\infty )$:
$$ H(x,y):=\sum_{x<p\le y}\frac{1}{p}
 - \log\log \theta(y^+) +\log\log \theta(x);
\ H(x):= \lim_{y\to\infty} H(x,y).
\eqno{(2.4)}
  $$

{\it Then}:
$$ {\bf (i)}\ -2.5 \le \liminf H(x)\, \sqrt{x}\, \log x; 
\ {\bf (ii)}\ \limsup H(x)\, \sqrt{x}\, \log x \le -1.5.
\eqno{(2.5)}
  $$

\vspace*{\fill}
\clearpage

$\triangleright$ \ First let us note that the Taylor formula implies:
$$ \log\log\theta(x) - \log\log x =\underline{\frac{\Delta(x)}{x\log x}}
- \frac{\log \hat{x} +1}{2\hat{x}^2\log^2\hat{x}} \Delta^2(x),
\eqno{(2.6)}
  $$

where $\hat{x}$ is a certain number in between of $x$ and $\theta(x)$.

\vspace{2ex}

On the other hand, integrating twice by parts, one obtains (cf. also (2.3)):
$$ \sum_{x<p\le y} \frac{1}{p} - \log\log y +\log\log x
=\int\limits_x^{y^+} \frac{d\Delta(t)}{t\log t}
=\frac{\Delta(t)}{t\log t}\Biggr|_x^{y^+}
- \int\limits_x^{y^+} \Delta(x) \left(\frac{1}{t\log t}\right)^{\prime} \, dt
$$
$$ =\frac{\Delta(y^+)}{y\log y}\ \ \underline{ -\ \frac{\Delta(x)}{x\log x}}
-\Phi(t)\left(\frac{1}{t\log t}\right)^{\prime}\Biggr|_x^y
 +\int\limits_x^y \Phi(t)\left(\frac{1}{t\log t}\right)^{\prime\prime}dt.
\eqno{(2.7)}
  $$

Making here $y$ tend to $\infty$, one comes to
$$ \lim_{y\to\infty}
 \left(\sum_{x<p\le y}\frac{1}{p} - \log\log y +\log\log x\right)
=\underline{ -\frac{\Delta(x)}{x\log x}} + D(x)+E(x)+F(x),
\eqno{(2.8)}
  $$
where the notations are used:
$$ D(x):=-\frac{\Phi(x)(\log x +1)}{x^2\log^2x},
 \ E(x):=\int\limits_x^{+\infty}\frac{\Phi(t)}{t^3\log t}
\left(2+\frac{3}{\log t}+\frac{2}{\log^2t}\right)\,dt,
$$
$$ F(x):=-\frac{\log\tilde{x}+1}{2\tilde{x}^2\log^2\tilde{x}}\, \Delta^2(x).
\eqno{(2.9)}
  $$

Summing this with (2.6) (the terms, involving $\Delta(x)$, underlined
in (2.6), (2.8), {\it mutually reduce}), and taking into account
 (2.2){\bf (iii)} and (2.3), one comes to: $ H(x)=D(x)+E(x)+O(\log^3x/x)$.
But by virtue of (2.3) one has:
$$\frac{17}{30}  \le \liminf D(x)\,\sqrt{x}\,\log x \le \limsup D(x)\,\sqrt{x}\,\log x \le\frac{23}{30};\qq\
$$
$$ -\frac{92}{30} \le \liminf E(x)\,\sqrt{x}\,\log x \le \limsup E(x)\,\sqrt{x}\,\log x \le - \frac{68}{30},
\eqno{(2.10)}
  $$
and thus $-2.5-\ve<H(x)\,\sqrt{x}\,\log x < -1.5+\ve$, for all $x>X_{\ve}$,
which coincides with (2.5) $\Box.$

\vspace*{\fill}
\clearpage

{\bf Remark 2.}\ It is important to emphasize that {\it without any
apriory estimates for} $\Delta(x)$, i. e. unconditionally,
from (2.8), (2.9) the inequality follows:
$$  H(x)\le D(x) + E(x), \q \forall x>3,
\eqno{(2.11)}
  $$
 because the quantity $F(x)$ is always non-positive.

Also it's easy to check that if the function $|b(x)|$ in (2.3) 
would be bounded by some $\delta_0>0$ (instead of 0.1),
then  the boundaries in (2.5) would be $-2\pm5\delta_0$.

\vspace{2ex}

\noindent
$\triangleright$\ Proceeding to the proof of the Theorem itself, let us note that defining formula (1.4) immediately implies:
\vspace{-2ex}
$$ \gamma= \lim_{y\to\infty}\left( \sum_{p\le y} 
 \log\frac{p}{p-1} -\log\log\theta(y)  \right),
\eqno{(2.12)}
  $$
and taking into account definitions (1.4) and (2.4), one obtains {\bf unconditionally}:
$$ Q(x)+H(x)=\sum_{p\le x}\log\frac{p}{p-1} - \gamma 
+ \lim_{y\to\infty}\left(\sum_{x<p\le y} \frac{1}{p} -\log\log\theta(y) \right) 
$$
\vspace{-2ex}
$$ =  \lim_{y\to\infty}\left(\sum_{p\le y}  \log\frac{p}{p-1} -\log\log\theta(y) \right)  - \gamma 
 + \sum_{p>x} \left(\frac{1}{p} -\log\frac{p}{p-1}\right)
$$
$$ = \sum_{p>x}  \left( \log\left(1-\frac{1}{p}\right)+\frac{1}{p}\right)=
 -\sum_{p>x}\sum_{k=2}^{\infty} \frac{1}{kp^k} 
=O\left(\frac{1}{x}\right),
\eqno{(2.13)}
  $$
and thus (2.5){\bf (i)(ii)} imply (1.8){\bf(i)(ii)} (respectively) $\Box.$
\vspace{1ex}

This completes  the Theorem's proof.

\vspace{3ex}

{\bf 3.\q Some corollaries and conclusive remarks}

\vspace{2ex}

The Theorem allows to deduce some new conditions equivalent to {\bf RH}
in terms of the function $\Phi(x)$ and of the primitive of $(\psi(x)-x)^2$.

\vspace{2ex}

{\bf Corollary 1.}\ {\it In order {\bf RH} hold true it is necessary
and sufficient that at least one (and then all) of the three conditions
be fulfilled:}
$$ {\bf (i)}\q \Phi(x) = O(x^{1.5+\ve}; 
\q {\bf (ii)}\ \Phi(x)=O(x^{1.5});
$$
$$ {\bf (iii)}\ \int\limits_0^x (\psi(t)-t)^2dt = O(x^2).
\eqno{(3.1)} 
$$

\noindent
$\triangleright$ \ {\it Necessity.} The implication {\bf RH} 
$\Rightarrow$ {\bf (iii)} was established by H. Cram\'er (1921),

 cf. [5], th. 13.5. 
Now it remains to notice that in (3.1) {\bf (iii)}$\Rightarrow$ {\bf (ii)}
 $\Rightarrow$ {\bf (i)}\ $\Box.$

\vspace*{\fill}
\clearpage

\vspace{1ex}

\noindent
$\triangleright$ \ {\it Sufficiency.}\ Let (3.1){\bf (i)} be fulfilled,
 i. e. $|\Phi(x)|\le C(\ve) x^{1.5+\ve}, \forall \ x, \ve>0$. 
then from (2.13), (2.11) and (2.9) one obtains for all $x>X_{\ve}$:
$$  -\ Q(x)=  H(x) +O\left(\frac{1}{x}\right)\le D(x)+E(x)
+O\left(\frac{1}{x}\right)
 $$
$$\le C_{\ve}\left(  \frac{x^{1.5+\ve}}{x^2\log^2x} 
+   2\int\limits_x^{\infty} \frac{t^{1.5+\ve}}{t^3\log t}\,dt \right)
\le C_{\ve}  x^{-0.5+\ve},
\eqno{(3.2)} 
$$
whence $Q(x)\ge - C_{\ve} x^{-0.5+\ve}, \ x>X_{\ve}$; but this coincides with (1.5){\bf(ii)}, which in turn  (by virtue of Proposition 1)
 implies {\bf RH} $\Box.$

\vspace{2ex}

The reasonings in the proofs of the Lemma and the Theorem show 
that
$$\hbox{\it both conditions  (2.5){\bf(i)(ii)} 
 are {\bf (separately) equivalent} to {\bf RH}}.
 \eqno{(3.3)} 
$$

This allows to deduce the  {\bf RH}-criteria 
 in terms of the consequences:
$$ U_k:= \sum_{j>k} \left(\, \frac{1}{p_j} - \frac{1}{\theta_j} \right),
\hbox{  where } \theta_j:=\theta(p_j); 
\q V_k:= \sum_{j>k} \left|\, \frac{1}{p_j} - \frac{1}{\theta_j}\, \right|.
\eqno{(3.4)} 
$$

\vspace{1ex}

{\bf Corollary 2.} {\it {\bf RH} is valid if and only if at least one (and then all) of the following six conditions  is fulfilled for any
 $\ve>0, \ k>K_{\ve}$:} $$ \qq  {\bf (i)}\  U_k < k^{-0.5+\ve};
 \qq \ {\bf (ii)}\  U_k > -  k^{-0.5+\ve}; \
$$
$$
\q  {\bf (iii)}\  U_k \sqrt{p_k}\,\log p_k \, < -1.5 + \ve;
\qq {\bf (iv)} \  U_k \sqrt{p_k}\,\log p_k \, > -2.5 - \ve;
$$
$$\qq\q {\bf (v)}\   V_k < k^{-0.5+\ve};\qq
\ {\bf (vi)}\  \frac{V_k \, \sqrt{p_k}}{\log p_k}
 < \frac{1+\ve}{4\pi}. 
\eqno{(3.5)} 
$$

\vspace{2ex}

\noindent
$\triangleright$ \ First we note that (cp. (2.6)) for  certain  $\tau_j\in (\theta_{j-1}, \theta_{j})$ one has:
$$
\log\log\theta_{j}-\log\log\theta_{j-1}
=\frac{\log p_j}{\theta_{j}\log\theta_j} 
+ \frac{\log\tau_j+1}{2\tau_j^2\log^2\tau_j} \log^2p_j 
$$
$$= \frac{1}{\theta_j}-\frac{\log\theta_j-\log p_j}{\theta_{j}\log\theta_j}
+O\left(\frac{1}{p_j^2\log p_j}\right) =  \frac{1}{\theta_j}
+O\left(\frac{1}{p_j^2\log p_j}\right) .
\eqno{(3.6)} 
$$

Here we have also taken into account that 
 $\log\theta_j-\log p_j=(\theta_j - p_j)/ \tilde{x}^*_j, $ where
$x^*_j$ is some number in between of $p_j$ and $\theta_j$,
and $\theta_j\approx p_j$.

Therefore, for $H_j:=H(p_j)$ (cf. (2.4)) one has:

\vspace*{\fill}
\clearpage

$$ H_{j-1} - H_{j} = \frac{1}{p_j} - \frac{1}{\theta_j}
+O\left(\frac{1}{p_j^2\log p_j}\right).
\eqno{(3.7)} 
$$

Summing these relations from $j=k$ to infinity, one obtains  
$$ H_k= U_k + O(1/k)= -Q(x)+O(1/k), \qq x\in(p_{k-1}, p_k]. 
\eqno{(3.8)} 
$$
and thus (3.5){\bf (i)(ii)(iii)(iv)}  are equivalent (resp.) to (1.5){\bf(i)(ii)},
 (1.8){\bf(i)(ii)}, each of which in turn (cf. Sect. 1) $\iff{\bf RH}$.

At last, taking  into account that $p_j\approx j\log j, \ j\to\infty,$ 
one obtains by virtue of (2.2){\bf(iii)}  under  
assumption of {\bf RH}:
$$ V_k = \sum_{j=k}^{\infty} \left| \frac{\theta_j-p_j}{\theta_jp_j}\right|
\le \frac{1+\delta_k}{8\pi} \sum_{j=k}^{\infty} \frac{\log^2p_j}{p_j\sqrt{p_j}}
$$
$$
\approx \frac{1}{8\pi} \sum_{j=k}^{\infty} \frac{\sqrt{\log j}}{j\sqrt{j}}
\approx  \frac{\sqrt{\log k}}{4\pi \sqrt{k}}
\approx\frac{\log p_k}{\sqrt{p_k}}, \ \hbox{ where } \delta_k\to0.
\eqno{(3.9)} 
$$

Hence, {\bf RH}$\Rightarrow$ (3.5){\bf (vi)} $\Rightarrow$ (3.5) {\bf (v)}
$\Rightarrow$ (3.5){\bf (i)}$\Rightarrow$ ${\bf RH} \ \Box$

\vspace{2ex}

{\bf Remark 3.} If {\bf RH} holds true, then Theorem shows not only
the fast, but also {\bf quasi-monotonic} decrease of the remainder
 $Q(x)$ in (1.4),  in the sense that relations $y>ax, a>25/9, x>X_a,$ 
imply $Q(x)>Q(y) > 0$. 

 The {\bf oscilation properties} of the remainder $R(x)$ in the original Mertens formula (1.2), studied by H. Diamond and J. Pintz [7], may be easily derived from the Theorem  for the  case, 
when {\bf RH} is valid.

Indeed, according to (2.1) one has {\bf unconditionally}:
$$R(x)=Q(x) + \frac{\Delta(x)}{\tilde{x} \log\tilde{x}}, 
\hbox{ where}\ \tilde{x}\approx x.
\eqno{(3.10)}
$$

Therefore {\it assuming\ {\bf RH}},  which according to the Theorem
(cf. (1.8)) implies $Q(x)=O(1/\sqrt{x} \log x)$, one obtains that
{\it for any positive function} 
$\eta(x)\to+\infty$: 
$${\bf(i)}\q \Delta(x) = \Omega_{\pm}(\sqrt{x}\, \eta(x)) \iff \ 
{\bf(ii)}\q R(x) = \Omega_{\pm}\left(\frac{\eta(x)}{\sqrt{x} \, \log(x)}\right)\ 
\eqno{(3.11)}
$$

By virtue of the J.\,Littlewood result  (1914), (cf., e. g. [5], th. 6.20),
this holds  true for $\eta(x)=\log\log\log x$. 

The same arguments allow to conclude that the  estimate 
 $R(x)=O(\log x / \sqrt{x})$, cf. (1.3), 
{\it is not only {\bf necessary}, but
also {\bf sufficient} for {\bf RH}}.

\vspace*{\fill}

\clearpage

\vspace{2ex}

\centerline{\bf LIST OF REFERENCES }

\vspace{3ex}

\noindent
[1] Mertens F. {\it $\ddot{U}ber$ einige asymptotische Gesetze der Zahlentheorie}.

 J. Reine Angew.  Math.,  {\bf 77}, 1874, pp. 289 -- 338.

\vspace{1ex}

\noindent
[2]\   Robin G. \ {\it Grandes valeurs de la fonction somme 
des diviseurs  }

{\it  et hypoth\`ese de  Riemann.}

 J. Math. Pures Appl. V. 63 (1984), pp. 187 -- 213.

\vspace{1ex}

\noindent
[3]\   Nicolas J.-L. \  {\it Petites valeurs de la fonction d'Euler.} 

 Journal of Number Theory, vol. 17, no. 3, 1983, p. 375 -- 388.

\vspace{1ex}

\noindent
[4]\ Kalyabin G.\,A. {\it Remainder in the Modified Mertens Formula

and Ramanujan Inequality}. 

 {\ https://arxiv.org/abs/2201.02663},  2022.

\vspace{1ex}

\noindent
[5]\ Narkiewicz W. {\it The Development of Prime Number Theory}. 
 
{Springer-Verlag Berlin Heidelberg New York},  2000.

\vspace{1ex}

\noindent
[6]\ Ingham A.\,E. {\it Distribution of Prime Numbers}. 
 
{\ Cambridge at the University Press},  1932.

\vspace{1ex}
\noindent
[7]\ Diamond H.\,G., Pintz J. {\it Oscilations of Mertens' product formula.}

Journal de Th\'eorie de Nombres de Bordeaux. T. 21, No 3 (2009), 523 - 533.

\vspace{1ex}

\end{document}